\documentclass[leqno,12pt]{amsart} 
\usepackage{amssymb}
\theoremstyle{plain} 
\newtheorem{theorem}{\indent\sc Theorem}[section]
\newtheorem{lemma}[theorem]{\indent\sc Lemma}

\theoremstyle{definition} 

\newtheorem{remark}[theorem]{\indent\sc Remark}

\begin{document}

\font\eightrm=cmr8
\font\eightit=cmti8
\font\eighttt=cmtt8
\def\nft
{\hbox{$n$\hskip3pt$\equiv$\hskip4pt$5$\hskip4.4pt$($mod\hskip2pt$3)$}}
\def\bbP{\text{\bf P}}
\def\bbR{\text{\bf R}}
\def\rto{\mathbf{R}\hskip-.5pt^2}
\def\rtr{\text{\bf R}\hskip-.7pt^3}
\def\rn{{\mathbf{R}}^{\hskip-.6ptn}}
\def\rk{{\mathbf{R}}^{\hskip-.6ptk}}
\def\bbZ{\text{\bf Z}}
\def\hyp{\hskip.5pt\vbox
{\hbox{\vrule width3ptheight0.5ptdepth0pt}\vskip2.2pt}\hskip.5pt}
\def\er{r\hskip.3pt}
\def\df{d\hskip-.8ptf}
\def\dfh{d\hskip-.8pt\fh}
\def\tb{T\hskip-.3pt\bs}
\def\tab{{T\hskip.2pt^*\hskip-2.3pt\bs}}
\def\txm{{T\hskip-2pt_x\hskip-.9ptM}}
\def\tym{{T\hskip-2pt_y\hskip-.9ptM}}
\def\txn{{T\hskip-2pt_x\hskip-.9ptN}}
\def\txtm{{T\hskip-2pt_{x(t)}\hskip-.9ptM}}
\def\txom{{T\hskip-2pt_{x(0)}\hskip-.9ptM}}
\def\txtsm{{T\hskip-2pt_{x(t,s)}\hskip-.9ptM}}
\def\tm{{T\hskip-.3ptM}}
\def\tn{{T\hskip-.3ptN}}
\def\tam{{T^*\!M}}
\def\so{\mathfrak{so}\hh}
\def\gi{\mathfrak{g}} 
\def\hi{\mathfrak{h}}
\def\fh{f}
\def\hf{\varLambda}
\def\rc{c}
\def\rd{d}
\def\kri{\text{\rm Ker}\hskip1pt\ri}
\def\kr{\text{\rm Ker}\hskip1ptR}
\def\kw{\text{\rm Ker}\hskip2ptW}
\def\kb{\text{\rm Ker}\hskip1.5ptB}
\def\xc{\mathcal{X}_c}
\def\dz{\mathcal{D}}
\def\dzp{\dz^\perp}
\def\dzx{\dz\hskip-.3pt_x}
\def\dzy{\dz\hskip-.3pt_y}
\def\dzxp{\dz\hskip-.3pt_x{}\hskip-4.2pt^\perp}
\def\dzyp{\dz\hskip-.3pt_y{}\hskip-4.2pt^\perp}
\def\lz{\mathcal{L}}
\def\tlz{\tilde{\lz}}
\def\xe{\mathcal{E}}
\def\xz{\mathcal{X}}
\def\yz{\mathcal{Y}}
\def\zz{\mathcal{Z}}
\def\tim{\hskip1.5pt\widetilde{\hskip-1.5ptM\hskip-.5pt}\hskip.5pt}
\def\tig{\hskip.7pt\widetilde{\hskip-.7ptg\hskip-.4pt}\hskip.4pt}
\def\hm{\hskip1.9pt\widehat{\hskip-1.9ptM\hskip-.2pt}\hskip.2pt}
\def\hmt{\hskip1.9pt\widehat{\hskip-1.9ptM\hskip-.5pt}_t}
\def\hmz{\hskip1.9pt\widehat{\hskip-1.9ptM\hskip-.5pt}_0}
\def\hmp{\hskip1.9pt\widehat{\hskip-1.9ptM\hskip-.5pt}_p}
\def\hg{\hskip1.2pt\widehat{\hskip-1.2ptg\hskip-.4pt}\hskip.4pt}
\def\hri{\hskip.7pt\widehat{\hskip-.7pt\ri\hskip-.4pt}\hskip.4pt}
\def\hn{\hskip.7pt\widehat{\hskip-.7pt\nabla\hskip-.4pt}\hskip.4pt}
\def\nao{\hbox{$\nabla\!\!^{^{^{_{\!\!\circ}}}}$}}
\def\ro{\hbox{$R\hskip-4.5pt^{^{^{_{\circ}}}}$}{}}
\def\mppp{\hbox{$-$\hskip1pt$+$\hskip1pt$+$\hskip1pt$+$}}
\def\mpdp{\hbox{$-$\hskip1pt$+$\hskip1pt$\dots$\hskip1pt$+$}}
\def\mmpp{\hbox{$-$\hskip1pt$-$\hskip1pt$+$\hskip1pt$+$}}
\def\mmmp{\hbox{$-$\hskip1pt$-$\hskip1pt$-$\hskip1pt$+$}}
\def\pppp{\hbox{$+$\hskip1pt$+$\hskip1pt$+$\hskip1pt$+$}}
\def\mpmp{\hbox{$-$\hskip1pt$\pm$\hskip1pt$+$}}
\def\mpmpp{\hbox{$-$\hskip1pt$\pm$\hskip1pt$+$\hskip1pt$+$}}
\def\mmpmp{\hbox{$-$\hskip1pt$-$\hskip1pt$\pm$\hskip1pt$+$}}
\def\q{q}
\def\bq{\hat q}
\def\p{p}
\def\w{\vt^\perp}
\def\x{v}
\def\y{y}
\def\vp{\vt^\perp}
\def\vd{\vt\hh'}
\def\vdx{\vd{}\hskip-4.5pt_x}
\def\bz{b\hh}
\def\cy{{y}}
\def\rwo{\,\hs\text{\rm rank}\hskip2.7ptW\hskip-2.7pt=\hskip-1.2pt1}
\def\rwho{\,\hs\text{\rm rank}\hskip2.2ptW^h\hskip-2.2pt=\hskip-1pt1}
\def\rw{\,\hs\text{\rm rank}\hskip2.4ptW\hskip-1.5pt}
\def\im{\varPhi}
\def\js{J}
\def\ism{H}
\def\fe{F}
\def\fy{f}
\def\dfc{dF\hskip-2.3pt_\cy\hskip.4pt}
\def\dfct{dF\hskip-2.3pt_\cy(t)\hskip.4pt}
\def\dic{d\im\hskip-1.4pt_\cy\hskip.4pt}
\def\vl{\Lambda}
\def\qt{\mathcal{E}}
\def\tqt{\tilde{\qt}}
\def\vh{h}
\def\mv{V}
\def\vy{\mathcal{V}}
\def\xv{\mathcal{X}}
\def\yv{\mathcal{Y}}
\def\iv{\mathcal{I}}
\def\gkp{\Sigma}
\def\bs{\varSigma}
\def\hs{\hskip.7pt}
\def\hh{\hskip.4pt}
\def\nh{\hskip-.7pt}
\def\nnh{\hskip-1pt}
\def\hrz{^{\hskip.5pt\text{\rm hrz}}}
\def\vrt{^{\hskip.2pt\text{\rm vrt}}}
\def\th{\varTheta}
\def\zh{\zeta}
\def\vg{\varGamma}
\def\my{\mu}
\def\ny{\nu}
\def\gy{\lambda}
\def\gm{\gamma}
\def\gp{\mathrm{G}}
\def\hp{\mathrm{H}}
\def\kp{\mathrm{K}}
\def\Gm{\Gamma}
\def\Lm{\Lambda}
\def\Dt{\Delta}
\def\sj{\sigma}
\def\lg{\langle}
\def\rg{\rangle}
\def\lr{\lg\,,\rg}
\def\uv{\underline{v\hskip-.8pt}\hskip.8pt}
\def\uvp{\underline{v\hh'\hskip-.8pt}\hskip.8pt}
\def\uw{\underline{w\hskip-.8pt}\hskip.8pt}
\def\uxs{\underline{x_s\hskip-.8pt}\hskip.8pt}
\def\vs{vector space}
\def\rvs{real vector space}
\def\vf{vector field}
\def\tf{tensor field}
\def\tvn{the vertical distribution}
\def\dn{distribution}
\def\od{Ol\-szak distribution}
\def\pt{point}
\def\tc{tor\-sion\-free connection}
\def\ea{equi\-af\-fine}
\def\rt{Ricci tensor}
\def\pde{partial differential equation}
\def\pf{projectively flat}
\def\pfs{projectively flat surface}
\def\pfc{projectively flat connection}
\def\pftc{projectively flat tor\-sion\-free connection}
\def\su{surface}
\def\sco{simply connected}
\def\psr{pseu\-\hbox{do\hs-}Riem\-ann\-i\-an}
\def\inv{-in\-var\-i\-ant}
\def\trinv{trans\-la\-tion\inv}
\def\feo{dif\-feo\-mor\-phism}
\def\feic{dif\-feo\-mor\-phic}
\def\feicly{dif\-feo\-mor\-phi\-cal\-ly}
\def\Feicly{Dif\-feo\-mor\-phi\-cal\-ly}
\def\diml{-di\-men\-sion\-al}
\def\prl{-par\-al\-lel}
\def\skc{skew-sym\-met\-ric}
\def\sky{skew-sym\-me\-try}
\def\Sky{Skew-sym\-me\-try}
\def\dbly{-dif\-fer\-en\-ti\-a\-bly}
\def\cs{con\-for\-mal\-ly symmetric}
\def\cf{con\-for\-mal\-ly flat}
\def\ls{locally symmetric}
\def\ecs{essentially con\-for\-mal\-ly symmetric}
\def\rr{Ric\-ci-re\-cur\-rent}
\def\kf{Killing field}
\def\om{\omega}
\def\vol{\varOmega}
\def\dv{\delta}
\def\ve{\varepsilon}
\def\zt{\zeta}
\def\kx{\kappa}
\def\mf{manifold}
\def\mfd{-man\-i\-fold}
\def\bmf{base manifold}
\def\bd{bundle}
\def\tbd{tangent bundle}
\def\ctb{cotangent bundle}
\def\bp{bundle projection}
\def\prc{pseu\-\hbox{do\hs-}Riem\-ann\-i\-an metric}
\def\prd{pseu\-\hbox{do\hs-}Riem\-ann\-i\-an manifold}
\def\Prd{pseu\-\hbox{do\hs-}Riem\-ann\-i\-an manifold}
\def\npd{null parallel distribution}
\def\pj{-pro\-ject\-a\-ble}
\def\pd{-pro\-ject\-ed}
\def\lcc{Le\-vi-Ci\-vi\-ta connection}
\def\vb{vector bundle}
\def\vbm{vec\-tor-bun\-dle morphism}
\def\kerd{\text{\rm Ker}\hskip2.7ptd}
\def\ro{\rho}
\def\sy{\sigma}
\def\ts{total space}
\def\pmb{\pi}
\def\ri{{\rho}}

\renewcommand{\thepart}{\Roman{part}}

\title[Con\-for\-mal\-ly symmetric manifolds]{The local structure of 
con\-for\-mal\-ly\\
symmetric manifolds}
\author[A. Derdzinski]{Andrzej Derdzinski}

\author[W. Roter]{Witold Roter}

%
\thanks{}

\address{
Department of Mathematics \endgraf
The Ohio State University \endgraf
Columbus, OH 43210 \endgraf
USA
}
\email{andrzej@math.ohio-state.edu}

\address{
Institute of Mathematics and Computer Science \endgraf
Wroc\l aw University of Technology \endgraf
Wy\-brze\-\.ze Wys\-pia\'n\-skiego 27, 50-370 Wroc\l aw \endgraf
Poland
}
\email{roter@im.pwr.wroc.pl}

\begin{abstract}This is a final step in a local classification of \prd s 
with parallel Weyl tensor that are not con\-for\-mal\-ly flat or locally 
symmetric.
\end{abstract}

\maketitle

\voffset=-18pt\hoffset=-4pt 

\setcounter{theorem}{0}
\renewcommand{\thetheorem}{\thesection.\arabic{theorem}}
\section*{Introduction}
The present paper provides a finishing touch in a local classification of 
\ecs\ \prc s.

A \prd\ of dimension $\,n\ge4\,$ is called {\it essentially con\-for\-mal\-ly 
symmetric} if it is {\it con\-for\-mal\-ly symmetric\/} \cite{chaki-gupta} 
(in the sense that its Weyl conformal tensor is parallel) without being 
con\-for\-mal\-ly flat or locally symmetric.

The metric of an \ecs\ \mf\ is always indefinite 
\cite[Theorem~2]{derdzinski-roter-77}. {\it Compact\/} \ecs\ manifolds are 
known to exist in all dimensions $\,n\ge5\,$ with $\,\nft$, where they 
represent all indefinite metric signatures \cite{derdzinski-roter}, while 
examples of \ecs\ \prc s on open manifolds of all dimensions $\,n\ge4\,$ were 
first constructed in \cite{roter}.

On every \cs\ \mf\ there is a naturally distinguished parallel distribution 
$\,\dz$, of some dimension $\,\rd$, which we call the {\it Ol\-szak 
distribution}. As shown by Ol\-szak \cite{olszak}, for an \ecs\ \mf\ 
$\,\rd\in\{1,2\}$.

In \cite{derdzinski-roter-07} we described the local structure of all \cs\ 
\mf s with $\,\rd=2$. See also Section~\ref{tcdt}. This paper establishes an 
analogous result (Theorem~\ref{clsdo}) for the case $\,\rd=1$.

In both cases, some of the metrics in question are locally symmetric. In 
Remark~\ref{njecs} we explain why a similar classification result cannot be 
valid just for {\it essentially\/} con\-for\-mal\-ly symmetric \mf s.

Essentially con\-for\-mal\-ly symmetric manifolds with $\,\rd=1\,$ are all 
{\it Ric\-ci-re\-cur\-rent}, in the sense that, for every tangent \vf\ $\,v$, 
the Ricci tensor $\,\ri\,$ and the covariant derivative $\,\nabla_{\!v}\ri\,$ 
are linearly dependent at each \pt. The local structure of \ecs\ \rr\ \mf s at 
points with $\,\ri\otimes\nnh\nabla\nnh\ri\ne0\,$ has already been 
determined by the second author \cite{roter}. Our new contribution settles the 
one case still left open in the local classification problem, namely, 
that of \ecs\ \mf s with $\,\rd=1\,$ at points where 
$\,\ri\otimes\nnh\nabla\nnh\ri=0$.

The literature dealing with \cs\ \mf s includes, among others, 
\cite{deszcz,deszcz-hotlos,hotlos,rong,sharma,simon} and the papers cited 
above. A local classification of {\it homogeneous} \ecs\ \mf s can be found in 
\cite{derdzinski-78}.

\section{Preliminaries}\label{prel}
Throughout this paper, all manifolds and bundles, along with sections and 
connections, are assumed to be of class $\,C^\infty\nnh$. A manifold is, 
by definition, connected. Unless stated otherwise, a mapping is always a 
$\,C^\infty$ mapping betweeen manifolds.

Given a connection $\,\nabla\,$ in a \vb\ $\,\xe\hs$ over a \mf\ $\,M$, 
a section $\,\psi\,$ of $\,\xe$, and \vf s $\,u,v\,$ tangent to $\,M$, 
we use the sign convention
\begin{equation}\label{cur}
R(u,v)\psi\hskip7pt=\hskip7pt\nabla_{\!v}\nabla_{\!u}\psi\,
-\,\nabla_{\!u}\nabla_{\!v}\psi\,+\,\nabla_{[u,v]}\psi
\end{equation}
for the curvature tensor $\,R=R^\nabla\nnh$.

The \lcc\ of a given \prd\ $\,(M,g)\,$ is always denoted by $\,\nabla$. 
We also use the symbol $\,\nabla\,$ for connections induced by $\,\nabla,$ 
in various $\,\nabla$\prl\ subbundles of $\,\tm\,$ and their quotients.

The Schouten tensor $\,\sigma\,$ and Weyl conformal tensor $\,W\hs$ of a \prd\ 
$\,(M,g)\,$ of dimension $\,n\ge4\,$ are given by 
$\,\sigma=\ri\,-\hs(2n-2)^{-1}\,\text{\rm s}\hskip1.2ptg$, with $\,\ri$ 
denoting the Ric\-ci tensor, 
$\,\hs\text{\rm s}\hs=\hs\text{\rm tr}_g\hh\ri\hs\,$ standing for the 
scalar curvature, and
\begin{equation}\label{wer}
W\,=\,\hs R\,-\,(n-2)^{-1}\hs g\wedge\hh\sigma\hs.
\end{equation}
Here $\,\wedge\,$ is the exterior multiplication of $\,1$-forms valued in 
$\,1$-forms, which uses the ordinary $\,\wedge\,$ as the val\-ue\-wise 
multiplication; thus, $\,g\wedge\hh\sigma\,$ is a $\,2$-form valued in 
$\,2$-forms.

Let $\,(t,s)\mapsto x(s,t)\,$ be a fixed {\it variation of curves} in a \prd\ 
$\,(M,g)$, that is, an $\,M$-val\-ued $\,C^\infty$ mapping from a rectangle 
(product of intervals) in the $\,ts\hs$-plane. By a {\it vector field\/ $\,w\,$ 
along the variation} we mean, as usual, a section of the pull\-back of 
$\,\tm\,$ to the rectangle (so that $\,w(t,s)\in\txtsm$). Examples are $\,x_s$ 
and $\,x_t$, which assign to $\,(t,s)\,$ the velocity of the curve 
$\,t\mapsto x(t,s)\,$ or $\,s\mapsto x(t,s)\,$ at $\,s\,$ or $\,t$. Further 
examples are provided by restrictions to the variation of vector fields on 
$\,M$. The partial covariant derivatives of a vector field $\,w\,$ along the 
variation are the vector fields $\,w_t,\hs w_s$ along the variation, obtained 
by differentiating $\,w\,$ covariantly along the curves $\,t\mapsto x(t,s)\,$ 
or $\,s\mapsto x(t,s)$. Skipping parentheses, we write $\,w_{ts},\hs w_{stt}$, 
etc., rather than $\,(w_t){}_s,\hs ((w_s){}_t){}_t$ for higher-or\-der 
derivatives, as well as $\,x_{ss},\hs x_{st}$ instead of 
$\,(x_s){}_s,\hs (x_s){}_t$. One always has $\,w_{ts}=w_{st}+R(x_t,x_s)\hh w$, 
cf.\ \cite[formula (5.29) on p.\ 460]{dillen-verstraelen}, and, since the 
\lcc\ $\,\nabla\,$ is tor\-sion\-free, $\,x_{st}=x_{ts}$. Thus, whenever 
$\,(t,s)\mapsto x(s,t)\,$ is a variation of curves in $\,M$,
\begin{equation}\label{xts}
x_{tss}\,\,=\,\,x_{sst}\,+\,\hs R(x_t,x_s)x_s\hs.
\end{equation}

\section{The Ol\-szak distribution}\label{oldi}
The {\it Ol\-szak distribution} of a \cs\ \mf\ $\,(M,g)\,$ is the parallel 
subbundle $\,\dz\,$ of $\,\tm$, the sections of which are the vector 
fields $\,u\,$ with the property that $\,\xi\wedge\hs\varOmega\hs=0\,$ for all 
vector fields $\,v,v\hh'$ and for the differential forms 
$\,\xi=g(u,\,\cdot\,)\,$ and 
$\,\hs\varOmega\hs=W(v,v\hh'\nnh,\,\cdot\,,\,\cdot\,)$. The distribution 
$\,\dz\,$ was introduced, in a more general situation, by Ol\-szak 
\cite{olszak}, who also proved the following lemma.
\begin{lemma}\label{oldis}The following conclusions hold for the dimension 
$\,\rd\,$ of the Ol\-szak distribution $\,\dz\,$ in any \cs\ manifold\/ 
$\,(M,g)\,$ with $\,\dim M=n\ge4$.
\begin{enumerate}
  \def\theenumi{{\rm\roman{enumi}}}
\item $\rd\in\{0,1,2,n\}$, and\/ $\,\rd=n\,$ if and only if\/ $\,(M,g)\,$ is 
con\-for\-mal\-ly flat.
\item $\rd\in\{1,2\}\,$ if\/ $\,(M,g)\,$ is \ecs.
\item $\rd=2\,$ if and only if\/ $\rwo$, in the sense that\/ $\,W\nh$, as an 
operator acting on exterior $\,2$-forms, has rank\/ $\,1\,$ at each \pt.
\item If\/ $\,\rd=2$, the distribution $\,\dz\,$ is spanned by all \vf s of 
the form $\,W(u,v)v'$ for arbitrary \vf s $\,u,v,v'$ on $\,M$.
\end{enumerate}
\end{lemma}
\begin{proof}See Appendix I.
\end{proof}
In the next lemma, parts (a) and (d) are due to Ol\-szak 
\cite[2$^{\hs\text{\rm o}}$ and 3$^{\hs\text{\rm o}}$ on p.~214]{olszak}.
\begin{lemma}\label{dontw}If\/ $\,\rd\in\{1,2\}$, where $\,\rd\,$ is 
the dimension of the Ol\-szak distribution\/ $\,\dz$ of a given \cs\ 
manifold\/ $\,(M,g)\,$ with $\,\dim M=n\ge4$, then
\begin{enumerate}
  \def\theenumi{{\rm\alph{enumi}}}
\item $\dz\,$ is a \npd,
\item at any $\,x\in M\,$ the space $\,\dz_x$ contains the image of the 
Ric\-ci tensor $\,\ri_x$ treated, with the aid of\/ $\,g_x$, as an 
endomorphism of\/ $\,\txm\nh$,
\item the scalar curvature is identically zero and\/ 
$\,R\hs\,=\,\hs W\hs+\,(n-2)^{-1}\hs g\wedge\ri$,
\item $W(u,\,\cdot\,,\,\cdot\,,\,\cdot\,)\hs\,=\,\hs0\hs\,$ whenever $\,u\,$ 
is a section of\/ $\,\dz$,
\item $R(v,v\hh'\nnh,\,\cdot\,,\,\cdot\,)
=\hs W(v,v\hh'\nnh,\,\cdot\,,\,\cdot\,)\hs\,=\,\hs0\hs\,$ for any sections 
$\,v\,$ and\/ $\,v\hh'$ of\/ $\,\dzp\nnh$,
\item of the connections in $\,\dz\hs$ and\/ $\,\qt=\dzp\nnh/\dz$, induced by 
the \lcc\ of\/ $\,g$, the latter is always flat, and the former is flat if\/ 
$\,\rd=1$.
\end{enumerate}
\end{lemma}
\begin{proof}Assertion (e) for $\,W\hs$ is immediate from the definition of 
$\,\dz$. \hbox{Namely, at any} point $\,x\in M$, every $\,2$-form 
$\,\hs\varOmega_x$ in the image of $\,W_{\nh x}$ (for $\,W_{\nh x}$ acting on 
\hbox{$\,2$-forms at $\,x$)} is $\,\wedge$-di\-vis\-i\-ble by 
$\,\xi=g_x(u,\,\cdot\,)\,$ for each 
$\,u\in\dzx\nh\smallsetminus\{0\}$, and so 
\hbox{$\,\varOmega_x(v,v\hh'\hh)=0\,$ if $\,v,v\hh'\nnh\in\nh\dzxp\nnh$.}

We now proceed to prove (a), (b), (c) and (d).

First, let $\,\rd=2$. By Lemma~\ref{oldis}(iii), this amounts to the condition 
$\rwo$, so that (a), (b) and (c) follow from Lemma~\ref{oldis}(iv) combined 
with \cite[Lemma 17.1(ii) and Lemma 17.2]{derdzinski-roter-07}. Also, for a 
nonzero $\,2$-form $\,\hs\varOmega_x$ chosen as in the last paragraph, 
$\,\dzx$ is the image of $\,\hs\varOmega_x$, that is, $\,\hs\varOmega_x$ 
equals the exterior product of two vectors in $\,\dzx$ (treated as 
$\,1$-forms, with the aid of $\,g_x$). Now (d) follows since, by (a), 
$\,\hs\varOmega_x(u_x,\,\cdot\,)=0$ if $\,u\,$ is a section of $\,\dz$.

Next, suppose that $\,\rd=1$. Replacing $\,M\,$ by a neighborhood of any given 
point, we may assume that $\,\dz\,$ is spanned by a vector field $\,u$. If 
$\,u\,$ were not null, we would have $\,W(u,v,u,v\hh'\hh)=0\,$ for any 
sections $\,v,v\hh'$ of $\,\dzp\nnh$, as one sees contracting the 
twice-co\-var\-i\-ant tensor field $\,W(\,\cdot\,,v,\,\cdot\,,v\hh'\hh)=0$, at 
any point $\,x$, in an orthogonal basis containing the vector $\,u_x$. (We 
have already established (e) for $\,W\nh$.) Combined with (e) for $\,W\hs$ 
and the symmetries of $\,W\nh$, the relation $\,W(u,v,u,v\hh'\hh)=0\,$ for 
$\,v,v\hh'$ in $\,\dzp$ would then give $\,W\nh=0$, contrary to the assumption 
that $\,\rd=1$. Thus, $\,u\,$ is null, which yields (a). Now
\begin{equation}\label{uuo}
\text{\rm we\ choose,\ locally,\ a\ null\ vector\ field\ $\,u\hh'$\ 
with\ $\,g(u,u\hh'\hh)=1$.}
\end{equation}
For any section $\,v\,$ of $\,\dzp$ one sees that 
$\,W(u,\,\cdot\,,u\hh'\nnh,v)=0\,$ by contracting the tensor field 
$\,W(\,\cdot\,,\,\cdot\,,\,\cdot\,,v)=0\,$ in the first and third arguments, 
at any point $\,x$, in
\begin{equation}\label{wtt}
\text{\rm a\ basis\ of\ $\,\txm\,$\ formed\ by\ $\,u_x,u_x'$\ and\ $\,n-2\,$\ 
vectors\ orthogonal\ to\ them,}
\end{equation}
and using (e) for $\,W\nh$, along with the inclusion $\,\dz\subset\dzp\nnh$, 
cf.\ (a). Since $\,u\hh'$ and $\,\dzp$ span $\,\tm$, assertion (e) for $\,W\hs$ 
thus implies (d).

To prove (b) and (c) when $\,\rd=1$, we distinguish two cases: $\,(M,g)\,$ is 
either \ecs, or locally symmetric. For (c), it suffices to establish vanishing 
of the scalar curvature $\,\hs\text{\rm s}\hs\,$ (cf.\ (\ref{wer})). Now, in 
the former case, $\,\hs\text{\rm s}\hs=0\,$ according to 
\cite[Theorem~7]{derdzinski-roter-78}, while (b) follows since, as shown in 
\cite[Theorem 7 on p.\ 18]{derdzinski-roter-80}, for arbitrary vector fields 
$\,v,v\hh'$ and $\,v\hh''$ on an \ecs\ \prd, $\,\xi\wedge\hs\varOmega\hs=0$, 
where $\,\xi=\ri\hh(v,\,\cdot\,)\,$ and 
$\,\hs\varOmega\hs=W(v\hh'\nnh,v\hh''\nnh,\,\cdot\,,\,\cdot\,)$. In the case 
where $\,g\,$ is locally symmetric, (b) and (c) are established in Appendix II.

Assertion (e) for $\,R\,$ is now obvious from (e) for $\,W\hs$ and (c), since, 
by (b), $\,\ri(v,\,\cdot\,)=0\,$ for any section $\,v\,$ of $\,\dzp\nnh$. The 
claim about $\,\qt\hs$ in (f) is in turn immediate from (\ref{cur}) and (e) 
for $\,R$, which states that $\,R(w,w\hh'\hh)\hh v$, for arbitrary vector 
fields $\,w,w\hh'$ and any section $\,v\,$ of $\,\dzp$, is orthogonal to all 
sections of $\,\dzp$ (and hence must be a section of $\,\dz$). Finally, to 
prove (f) for $\,\dz$, with $\,\rd=1$, let us fix a section $\,u\,$ of 
$\,\dz$, a \vf\ $\,v$, and define a differential $\,2$-form $\,\zeta\,$ by 
$\hs\zeta(w,w\hh'\hh)=(n-2)\hs R(w,w\hh'\nnh,u,v)\hs$ for any vector fields 
$\,w,w\hh'\nnh$. By (c) and (e), 
$\,\zeta=g(u,\,\cdot\,)\wedge\ri(v,\,\cdot\,)$, as $\,\dz\subset\dzp$ (cf.\ 
(a)), and so $\,\ri(u,\,\cdot\,)=0\,$ in view of (b) and symmetry of $\,\ri$.
However, by (b), both $\,g(u,\,\cdot\,)\,$ and $\,\ri(v,\,\cdot\,)\,$ are 
sections of the subbundle of $\,\tam\,$ corresponding to $\,\dz\,$ under the 
bundle isomorphism $\,\tm\to\tam\,$ induced by $\,g$, so that $\,\zeta=0\hs\,$ 
since the distribution $\,\nh\dz\hs$ is one\diml.
\end{proof}

\section{The case $\,\rd=2$}\label{tcdt}
For more details of the construction described below, we refer the reader 
to \cite{derdzinski-roter-07}.

Let there be given a \su\ $\,\bs$, a \pftc\ $\,\hs{\rm D}\hs\,$ on $\,\bs$ 
with a $\,\hs\text{\rm D}\hs$\prl\ area form $\,\alpha$, an integer 
$\,n\ge4\hh$, a sign factor $\,\ve=\pm1$, a \rvs\ $\,\mv\hs$ of dimension 
$\,n-4\hh$, and a pseu\-\hbox{do\hs-}Euclid\-e\-an inner product $\,\lr\,$ 
on $\,\mv\nnh$.

We also assume the existence of a twice-con\-tra\-var\-i\-ant symmetric 
tensor field $\hs\hs T$ on $\,\bs\,$ with 
$\,\hs\text{\rm div}{}^{\hs\text{\rm D}}
(\text{\rm div}{}^{\hs\text{\rm D}}T)\,
+\,(\hh\ri^{\hs\text{\rm D}}\nnh,T\hh)=\hh\ve\,$ (in coordinates: 
$\,T^{jk}{}_{,\hh jk}+T^{jk}R_{jk}=\hh\ve$). Here 
$\,\hs\text{\rm div}{}^{\hs\text{\rm D}}$ denotes the 
$\,\hs\text{\rm D}\hs$-di\-ver\-gence, $\,\ri^{\hs\text{\rm D}}$ is the \rt\ 
of $\,\hs\text{\rm D}\hh$, and $\,(\hskip2pt,\hskip1pt)\,$ stands for the 
obvious pairing. Such $\,T\,$ always exists locally in $\,\bs$. In fact, 
according to \cite[Theorem~10.2(i)]{derdzinski-roter-07} combined with 
\cite[Lemma~11.2]{derdzinski-roter-07}, $\,T\,$ exists whenever $\,\bs\,$ is 
simply connected and noncompact.

For $\,T\,$ chosen as above, we define a twice-co\-var\-i\-ant symmetric 
tensor field $\,\tau\,$ on $\,\bs$, that is, a section of 
$\,[\tab]^{\odot2}\nnh$, by requiring $\,\tau\,$ to correspond to the section 
$\,T\,$ of $\,[\tb]^{\odot2}$ under the vec\-tor-bun\-dle isomorphism 
$\,\tb\to\tab\,$ which acts on \vf s $\,v\,$ by 
$\,v\mapsto\alpha(v,\,\cdot\,)$. In coordinates, 
$\,\tau_{jk}=\alpha_{jl}\hh\alpha_{km}T^{\hs lm}\nnh$.

Next, we denote by $\,h^{\text{\rm D}}$ the {\it Patter\-son\hs-\nh Walk\-er 
Riemann extension metric} \cite{patterson-walker} on the total space $\,\tab$, 
obtained  by requiring that all vertical and all 
$\,\hs\text{\rm D}\hs$-hor\-i\-zon\-tal vectors be 
$\,h^{\text{\rm D}}\nnh$-null, while 
$\,h_x^{\text{\rm D}}(\zeta,w)=\zeta(d\pmb_xw)\,$ for $\,x\in\tab$, any vector 
$\,w\in T_x\tab$, any vertical vector 
$\,\zeta\in\kerd\pmb_x=T_{\nnh\pmb(x)}^*\hskip-1pt\bs$, and the \bp\ 
$\,\pmb:\tab\to\bs$.

Finally, let $\,\gm\,$ and $\,\theta\,$ be the constant \prc\ on $\,\mv$ 
corresponding to the inner product $\,\lr$, and the function 
$\,\mv\nh\to\bbR\,$ with $\,\theta(v)=\lg v,v\rg$.

Our $\,\bs,\hs\text{\rm D}\hh,\alpha,n,\ve,\mv\nh,\lr\,$ now give rise to the 
\prd
\begin{equation}\label{hgt}
(\tab\,\times\,\mv,\,\,h^{\text{\rm D}}\hskip-1.9pt-2\tau+\gm
-\theta\ri^{\hs\text{\rm D}})\,,
\end{equation}
of dimension $\,n$, with the metric 
$\,h^{\text{\rm D}}\hskip-2.5pt-2\tau+\gm-\theta\ri^{\hs\text{\rm D}}\nnh$, 
where the function $\,\theta\,$ and covariant tensor fields 
$\,\tau,\ri^{\hs\text{\rm D}}\nnh,h^{\text{\rm D}}\nnh,\gm\,$ on 
$\,\bs,\,\tab\,$ or $\,\mv\nh$ are identified with their pull\-backs to 
$\,\tab\,\times\,\mv\nnh$. (Thus, for instance, 
$\,h^{\text{\rm D}}\hskip-1.9pt-2\tau+\gm\,$ is a product metric.)

We have the following local classification result, in which $\,\rd\,$ stands 
for the dimension of Ol\-szak distribution $\,\dz$.
\begin{theorem}\label{clsdt}The \prd\/ {\rm(\ref{hgt})} obtained as above from 
any data $\,\bs,\hs\text{\rm D}\hh,\alpha,n,\ve,\mv\nh,\lr\,$ with the stated 
properties is \cs\ and has\/ $\,\rd=2$. Conversely, in any \cs\ \prd\ such 
that\/ $\rd=2$, every \pt\ has a connected neighborhood isometric to an open 
subset of a manifold\/ {\rm(\ref{hgt})} constructed above from some data 
$\,\bs,\hs\text{\rm D}\hh,\alpha,n,\ve,\mv\nh,\lr$.

The manifold\/ {\rm(\ref{hgt})} is never con\-for\-mal\-ly flat, and it 
is locally symmetric if and only if the \rt\ $\,\ri^{\hs\text{\rm D}}$ is 
$\,\hs\text{\rm D}\hs$\prl.
\end{theorem}
\begin{proof}See \cite[Section~22]{derdzinski-roter-07}. Note that, in view of 
Lemma~\ref{oldis}(iii), the condition $\rwo\,$ used in 
\cite{derdzinski-roter-07} is equivalent to $\rd=2$.
\end{proof}

The objects $\,\bs,\hs\text{\rm D}\hh,\alpha,n,\ve,\mv\nh,\lr\,$ are treated 
as parameters of the above construction, while $\,T\,$ is merely assumed to 
exist, even though the metric $\,g\,$ in (\ref{hgt}) clearly depends on 
$\,\tau\,$ (and hence on $\,T$). This is justified by the fact that, with 
fixed $\,\bs,\hs\text{\rm D}\hh,\alpha,n,\ve,\mv\nh,\lr$, the metrics 
corresponding to two choices of $\,T\,$ are, locally, isometric to each 
other, cf.\ \cite[Remark~22.1]{derdzinski-roter-07}.

The metric signature of (\ref{hgt}) is clearly given by 
$\,-\hs-\hs\ldots\hs+\hs+\hh$, with the dots standing for the sign pattern of 
$\,\lr$.

\section{The case $\,\rd=1$}\label{tcdo}
Let there be given an open interval $\,I\nh$, a $\,C^\infty$ function 
$\,\fh:I\to\bbR\hs$, an integer $\,n\ge4$, a \rvs\ $\,\mv\hs$ of dimension 
$\,n-2\,$ with a pseu\-\hbox{do\hs-}Euclid\-e\-an inner product $\,\lr$, and a 
nonzero traceless linear operator $\,A:\mv\to\mv\nh$, self-ad\-joint relative 
to $\,\lr$. As in \cite{roter}, we then define an $\,n$\diml\ \prd\ 
\begin{equation}\label{rcr}
(I\times\bbR\times\mv,\,\,\kx\,dt^2\hs+\,dt\,ds\,+\,\gm)\,,
\end{equation}
where products of differentials represent symmetric products, $\,t,s\,$ denote 
the Cartesian coordinates on the $\,I\times\bbR\,$ factor, $\,\gm\,$ stands for 
the pull\-back to $\,I\times\bbR\times\mv\hs$ of the flat \prc\ on $\,\mv$ that 
corresponds to the inner product $\,\lr$, and the function 
$\,\kx:I\times\bbR\times\mv\nh\to\bbR\,$ is given by 
$\,\kx(t,s,\psi)=\fh(t)\hh\lg \psi,\psi\rg+\lg A\psi,\psi\rg$.

The manifolds (\ref{rcr}) are characterized by the following local 
classification result, analogous to Theorem~\ref{clsdt}. As before, 
$\,\rd\,$ is the dimension of the Ol\-szak distribution.
\begin{theorem}\label{clsdo}For any $\,I\nh,\fh,n,\mv\nh,\lr,A\,$ as above, 
the \prd\/ {\rm(\ref{rcr})} is \cs\ and has\/ $\,\rd=1$. Conversely, in 
any \cs\ \prd\ such that\/ $\rd=1$, every \pt\ has a connected neighborhood 
isometric to an open subset of a manifold\/ {\rm(\ref{rcr})} constructed from 
some such $\,I\nh,\fh,n,\mv\nh,\lr,A$.

The manifold\/ {\rm(\ref{rcr})} is never con\-for\-mal\-ly flat, and it 
is locally symmetric if and only if\/ $\,\fh\,$ is constant.
\end{theorem}
A proof of Theorem~\ref{clsdo} is given at the end of the next section.

Obviously, the metric $\,\,\kx\,dt^2\hs+\,dt\,ds\,+\,\gm\,\,$ in (\ref{rcr}) 
has the sign pattern $\,\,-\hs\ldots\hs+\hh$, where the dots stand for the 
sign pattern of $\,\lr$.
\begin{remark}\label{njecs}A classification result of the same format as 
Theorem~\ref{clsdo} cannot be true just for {\it essentially\/} 
con\-for\-mal\-ly symmetric \mf s with $\,\rd=1$. Namely, such \mf s do not 
satisfy a principle of unique continuation: formula (\ref{rcr}) with 
$\,\fh$ which is nonconstant on $\,I\nh$, but constant on some nonempty 
open subinterval $\,I\hs'$ of $\,I\nh$, defines an \ecs\ \mf\ with a locally 
symmetric open submanifold $\,\,U=I\hs'\nnh\times\bbR\times\mv\nnh$. At 
points of $\,\,U\nh$, the local structure of (\ref{rcr}) does not, therefore, 
arise from a construction that, locally, produces all \ecs\ \mf s and nothing 
else.

As explained in \cite[Section~24]{derdzinski-roter-07}, an analogous 
situation arises when $\,\rd=2$.
\end{remark}

\section{Proof of Theorem~\ref{clsdo}}\label{potf}
The following assumptions will be used in Lemma~\ref{fstar}.
\begin{enumerate}
  \def\theenumi{{\rm\alph{enumi}}}
\item $(M,g)\,$ is a \cs\ \mf\ of dimension $\,n\ge4\,$ and $\,y\in M$.
\item The \od\ $\,\dz\,$ of $\,(M,g)\,$ is one\diml.
\item $u\,$ is a global parallel vector field spanning $\,\dz$. 
\item $t:M\to\bbR\,$ is a $\,C^\infty$ function with $\,g(u,\,\cdot\,)=dt\,$ 
and $\,t(y)=0$.
\item $\dim\mv\nh=n-2\,$ for the space $\,\mv\hs$ of all parallel sections of 
$\,\qt=\dzp\nnh/\dz$.
\item $\ri=(2-n)\fh(t)\,dt\otimes dt\,$ for some $\,C^\infty$ function 
$\,\fh:I\hh'\nh\to\bbR\,$ on an open interval $\,I\hh'\nnh$, where $\,\ri\,$ 
is the Ric\-ci tensor and $\,\fh(t)\,$ denotes the composite 
$\,\fh\nnh\circ\hh t$.
\end{enumerate}
For local considerations, only (a) and (b) are essential. In fact, condition 
(e) (in which `parallel' refers to the connection in $\,\qt\,$ induced by the 
\lcc\ of $\,g$), as well (c) and (d) for some $\,u\,$ and $\,t$, follow from 
(a) -- (b) if $\,M\,$ is simply connected. See Lemma~\ref{dontw}(f). On the 
other hand, (c) -- (d), Lemma~\ref{dontw}(b) and symmetry of $\,\ri\,$ give 
$\,\nabla dt=0\,$ and $\,\ri=\chi\,dt\otimes dt\,$ for some function 
$\,\chi:M\to\bbR\hs$, so that 
$\,\nabla\nh\ri=d\chi\otimes\hs dt\otimes\hs dt$. However, $\,\nabla\nh\ri\,$ 
is totally symmetric (that is, $\,\ri\,$ satisfies the Co\-daz\-zi equation): 
our assumption $\,\nabla\hs W\nnh=0\,$ implies the condition 
$\,\hs\text{\rm div}\,W\nnh=0$, well known 
\cite[formula (5.29) on p.\ 460]{dillen-verstraelen} to be equivalent to 
the Co\-daz\-zi equation for the Schouten tensor $\,\sigma$, while 
$\,\sigma=\ri\,$ by Lemma~\ref{dontw}(c). Thus, $\,d\chi\,$ equals 
a function times $\,dt$, and so $\,\chi\,$ is, locally, a function of $\,t$, 
which (locally) yields (f).

For any section $\,v\,$ of $\,\dzp\nnh$, we denote by $\,\uv\,$ the image of 
$\,v\,$ under the quo\-tient-pro\-jec\-tion morphism 
$\,\dzp\nnh\to\,\qt=\dzp\nnh/\dz$.

The data required for the construction in Section~\ref{tcdo} consist of 
$\,I\nh,\fh,n,\mv\hs$ appearing in (a) -- (f), along with the 
pseu\-\hbox{do\hs-}Euclid\-e\-an inner product $\,\lr\,$ in $\,\mv\nh$, 
induced in an obvious way by $\,g\,$ (cf.\ Lemma~\ref{dontw}(f)), and 
$\,A:\mv\to\mv\hs$ characterized by 
$\,\lg A\psi,\psi\hh'\hh\rg
=W(u\hh'\nnh,v,v\hh'\nnh,u\hh'\hh)$, for $\,\psi,\psi\hh'\nh\in\mv\nh$, with a 
vector field $\,u\hh'$ and sections $\,v,v\hh'$ of $\,\dzp$ chosen, locally, 
so that $\,g(u,u\hh'\hh)=1$, $\,\psi=\uv\,$ and $\,\psi\hh'\nh=\uvp$. (The 
resulting bilinear form 
$\,(\psi,\psi\hh'\hh)\mapsto\lg A\psi,\psi\hh'\hh\rg\,$ on $\,\mv\hs$ is 
well-de\-fined, that is, unaffected by the choices of $\,u\hh'\nnh,v\,$ or 
$\,v\hh'\nnh$, as a consequence of Lemma~\ref{dontw}(d),\hs(e), while the 
function $\,W(u\hh'\nnh,v,v\hh'\nnh,u\hh'\hh)\,$ is in fact constant, by 
Lemma~\ref{dontw}(d), as ones sees differentiating it via the Leib\-niz rule 
and noting that, since $\,\uv\,$ and $\,\uvp$ are parallel, the covariant 
derivatives of $\,v\,$ and $\,v\hh'$ in the direction of any vector field are 
sections of $\,\dz$.) That $\,A\,$ is traceless and self-ad\-joint is 
immediate from the symmetries of $\,W\nh$. Finally, $\,A\ne0\,$ since, 
otherwise, $\,W\hs$ would vanish. (Namely, in view of 
Lemma~\ref{dontw}(d),\hs(e), $\,W\hs$ would yield $\,0\,$ when evaluated on any quadruple of vector fields, 
each of which is either $\,u\hh'$ or a section of $\,\dzp\nnh$.)

Under the assumptions (a) -- (f), with $\,\fh=\fh(t)$, we then have
\begin{equation}\label{ruv}
R(u\hh'\nnh,v)\hh v\hh'\,=\,\hs[\fh \hs g(v,v\hh'\hh)
+\lg A\uv,\uvp\hh\rg]\hs g(u\hh'\nnh,u)\hs u
\end{equation}
for any sections $\,v,v\hh'$ of $\,\dzp$ and any vector field $\,u\hh'\nnh$. 
In fact, $\,\ri\hh(v,\,\cdot\,)=\ri\hh(v\hh'\nnh,\,\cdot\,)=0$ from symmetry 
of $\,\ri\,$ and Lemma~\ref{dontw}(b), so that, by Lemma~\ref{dontw}(c), 
$\,R(u\hh'\nnh,v)\hh v\hh'\nh=W(u\hh'\nnh,v)\hh v\hh'\nh
-(n-2)^{-1}g(v,v\hh'\hh)\ri u\hh'\nnh$, where $\,\ri u\hh'$ denotes the unique 
vector field with $\,g(\ri u\hh'\nnh,\,\cdot\,)=\ri\hh(u\hh'\nnh,\,\cdot\,)$. 
Now (\ref{ruv}) follows: due to (d), (f) and the definition of $\,A$, 
both sides have the same $\,g$-in\-ner product with $\,u\hh'\nh$, and are 
orthogonal to $\,u^\perp\nh=\dzp$ (with $\,R(u\hh'\nnh,v)\hh v\hh'$ orthogonal 
to $\,\dzp$ in view of Lemma~\ref{dontw}(e)).

We fix an open subinterval $\,I\hs$ of $\,I\hh'\nnh$, containing $\,0$, and a 
null geodesic $\,I\ni t\mapsto x(t)\,$ in $\,M\,$ with $\,x(0)=y$, 
parametrized by the function $\,t\,$ (in the sense that the function $\,t\,$ 
restricted to the geodesic coincides with the geodesic parameter). Namely, 
since $\,\nabla dt=0$, the restriction of $\,t\,$ to any geodesic is an affine 
function of the parameter; thus, by (d), it suffices to prescribe the initial 
data formed by $\,x(0)=y\,$ and a null vector $\,\dot x(0)\in\tym\,$ with 
$\,g(\dot x(0),u_y)=1$.

As $\,g(\dot x(0),u_y)=1$, the plane $\,P\,$ in $\,\tym$, spanned 
by the null vectors $\,\dot x(0)\,$ and $\,u_y$ (cf.\ Lemma~\ref{dontw}(a)) is 
$\,g_y$-non\-de\-gen\-er\-ate, and so $\,\tym=P\oplus\widetilde\mv\nh$, for 
$\,\widetilde\mv\nh=P^\perp\nnh$. Let 
$\,\hs\text{\rm pr}:\tym\to\widetilde\mv\hs$ be the orthogonal projection. 
Since $\,\hs\text{\rm pr}\hs(\dzy)=\{0\}$, the restriction of 
$\,\hs\text{\rm pr}\hs\,$ to $\,\dzyp$ descends to the quotient 
$\,\qt_y=\dzyp\nnh/\dz_y$, producing an isomorphism 
$\,\qt_y\to\widetilde\mv\nh$, also denoted by $\,\hs\text{\rm pr}\hh$. 
Finally, for $\,\psi\in\mv\nh$, we let $\,t\mapsto \tilde \psi(t)\in\txtm\,$ 
be the parallel field with 
$\,\tilde \psi(0)=\hs\text{\rm pr}\hskip2.4pt\psi_y$, and set 
$\,\kx(t,s,\psi)=\fh(t)\hh\lg \psi,\psi\rg+\lg A\psi,\psi\rg$, as in 
Section~\ref{tcdo}.

The formula $\,F(t,s,\psi)=\exp_{\hs x(t)}(\tilde \psi(t)+su_{x(t)}/2)\,$ now 
defines a $\,C^\infty$ mapping $\,F$ from an open subset of 
$\,\rto\nh\times\mv\hs$ into $\,M$.
\begin{lemma}\label{fstar}Under the above hypotheses, 
$\,F\hh^*\nnh g=\kx\,dt^2\nh+\hs dt\hs ds+\vh$.
\end{lemma}
\begin{proof}The $\,F$-im\-ages $\,w,w\hh'\nnh,F_*\psi\,$ of the constant 
vector fields $\,(1,0,0),(0,1,0)$ and $\,(0,0,\psi)\,$ in 
$\,\rto\nh\times\mv\nh$, for $\,\psi\in\mv\nh$, are vector fields tangent to 
$\,M\,$ along $\,F$ (sections of $\,F\hh^*\tm$). Since $\,\dzp$ is parallel, 
its leaves are totally geodesic and, by Lemma~\ref{dontw}(e), the \lcc\ of 
$\,g\,$ induces on each leaf a flat tor\-sion\-free connection. Thus, 
$\,w\hh'$ and each $\,F_*\psi\,$ are parallel along each leaf of $\,\dzp\nnh$, as 
well as tangent to the leaf, and parallel along the geodesic 
$\,t\mapsto x(t)$. Therefore, $\,w\hh'\nh=u/2$, while the 
functions $\,g(w\hh'\nnh,F_*\psi)\,$ and $\,g(F_*\psi,F_*\psi\hh'\hs)$, for 
$\,\psi,\psi\hh'\nh\in\mv\nh$, are constant, and hence equal to their values at 
$\,y$, that is, $\,0\,$ and $\,\lg \psi,\psi\hh'\hs\rg$. It now remains to be 
shown that $\,g(w,w)=\kx\circ F$, $\,g(w,u/2)=1/2\,$ and 
$\,g(w,F_*\psi)=0$. To this end, we consider the variation 
$\,x(t,s)=F(t,sa,s\psi)\,$ of curves in $\,M$, with any fixed $\,a\in\bbR\,$ and 
$\,\psi\in\mv\nh$. Clearly, $\,w=x_t$ along the variation (notation of 
Section~\ref{prel}). Next, $\,x_{ts}=x_{st}$ is tangent to $\,\dzp\nnh$, since 
so is $\,x_s$, while $\,\dzp$ is parallel. Consequently, 
$\,[\hs g(x_t,u)\hh]_s=0$, as $\,u\,$ is parallel and tangent to $\,\dz$. 
Thus, $\,g(w,u)=g(x_t,u)=1$. (Note that $\,g(x_t,u)=1\,$ at $\,s=0$, due to 
(d), as the geodesic $\,t\mapsto x(t)\,$ is parametrized by the function 
$\,t$.) However, $\,x_{ss}=0\,$ and $\,x_s$ is tangent to $\,\dzp\nnh$, so 
that (\ref{xts}) and (\ref{ruv}) now give 
$\,x_{tss}=[\hs\fh g(x_s,x_s)+\lg A\uxs,\uxs\rg\hh]\hs u$, which is parallel 
in the $\,s\,$ direction, while $\,x_{ts}=x_{st}=0\,$ at $\,s=0$. Hence 
$\,x_{ts}=s\hs[\hs\fh g(x_s,x_s)+\lg A\uxs,\uxs\rg\hh]\hs u$, and so 
$\,g(x_{ts},x_{ts})=0\,$ (cf.\ (c) above and Lemma~\ref{dontw}(a)). This 
further yields $\,[\hs g(x_t,x_t)]_{ss}/2=g(x_t,x_{tss})
=\fh g(x_s,x_s)+\lg A\uxs,\uxs\rg$. The last function is constant in the 
$\,s\,$ direction, while $\,g(x_t,x_t)=[\hs g(x_t,x_t)]_s=0\,$ at $\,s=0$, and 
so $\,g(w,w)=g(x_t,x_t)=s^2[\hs\fh g(x_s,x_s)+\lg A\uxs,\uxs\rg]=\kx$. 
Finally, being proportional to $\,u\,$ at each point, $\,x_{ts}$ is orthogonal 
to $\,\dzp\nnh$, and hence to $\,F_*\psi$, which imples that 
$\,[\hs g(x_t,F_*\psi)\hh]_s=0$, and, as $\,g(w,F_*\psi)=g(x_t,F_*\psi)=0\,$ 
at $\,s=0$, we get $\,g(w,F_*\psi)=0$ everywhere.
\end{proof}
We are now in a position to prove Theorem~\ref{clsdo}. First, (\ref{rcr}) is 
\cs\ and has $\,\rd=1$, as one can verify by a direct calculation, cf.\ 
\cite[Theorem 3]{roter}. Conversely, if conditions (a) and (b) above are 
satisfied, we may also assume (c) -- (f). (See the comment following (f).) Our 
assertion is now immediate from Lemma~\ref{fstar}.

\setcounter{section}{1}
\renewcommand{\thesection}{\Roman{section}}
\section*{Appendix I: Proof of Lemma~\ref{oldis}}\label{appo}
We prove Lemma~\ref{oldis} here, since Ol\-szak's paper \cite{olszak} may be 
difficult to obtain.

The condition $\,\rd=n\,$ is equivalent to con\-for\-mal flatness 
of $\,(M,g)$, since $\,n>2$ and so $\,\hs\varOmega\hs=0\,$ is the only 
$\,2$-form $\,\wedge$-di\-vis\-i\-ble by all nonzero $\,1$-forms $\,\xi$. At a 
fixed point $\,x$, the metric $\,g_x$ allows us to treat the Ric\-ci tensor 
$\,\ri_x$ and any $\,2$-form $\,\hs\varOmega_x$ as endomorphisms of 
$\,\txm\nh$, so that we may consider their images (which are subspaces of 
$\,\txm$). If $\,W\nh\ne0$, fixing a nonzero $\,2$-form $\,\hs\varOmega_x$ in 
the image of $\,W_{\nh x}$ acting on $\,2$-forms at $\,x\,$ we see that, for 
every $\,u\in\dzx$, our $\,\hs\varOmega_x$ is $\,\wedge$-di\-vis\-i\-ble by 
$\,\xi=g_x(u,\,\cdot\,)$, and so the image of $\,\hs\varOmega_x$ contains 
$\,\dzx$. Thus, $\,\rd\le2$, and (i) follows. (Being nonzero and decomposable, 
$\,\hs\varOmega_x$ has rank $\,2$.) As shown in \cite[Theorem 7 on p.\ 
18]{derdzinski-roter-80}, if $\,(M,g)$ is \ecs, the image of $\,\ri_x$ is a 
subspace of $\,\dzx$, so that (i) yields (ii), since $\,g\,$ in (ii) cannot be 
Ric\-ci-flat. Next, if $\,\rd=2$, the image of our $\,\hs\varOmega_x$ 
coincides with $\,\dzx$ (as $\,\hs\text{\rm rank}\,\varOmega_x=2$). Every 
$\,2$-form in the image of $\,W_{\nh x}$ thus is a multiple of 
$\,\hs\varOmega_x$, being the exterior product of two vectors in $\,\dzx$, 
identified, via $\,g_x$, with $\,1$-forms. Hence $\rwo$. Conversely, if 
$\rwo$, all nonzero $\,2$-forms $\,\hs\varOmega_x$ in the image of 
$\,W_{\nh x}$ are of rank $\,2$, as $\,W_{\nh x}$, being self-ad\-joint, is a 
multiple of $\,\hs\varOmega_x\otimes\hs\varOmega_x$, and so the Bianchi 
identity for $\,W\hs$ gives $\,\hs\varOmega_x\wedge\hs\varOmega_x=0$. All such 
$\,\hs\varOmega_x$ are therefore $\,\wedge$-di\-vis\-i\-ble by 
$\,\xi=g_x(u,\,\cdot\,)$, for every nonzero vector $\,u\,$ in the common 
$\,2$\diml\ image of such $\,\hs\varOmega_x$, which shows that $\,\rd=2$. 
Finally, (iv) follows if one chooses $\,\hs\varOmega_x\ne0\,$ equal to 
$\,W_{\nh x}(v,v\hh'\nnh,\,\cdot\,,\,\cdot\,)\,$ for some 
$\,v,v\hh'\in\txm$.

\setcounter{section}{2}
\renewcommand{\thesection}{\Roman{section}}
\section*{Appendix II: Lemma~\ref{dontw}{\rm(b),\hs(c)} in the locally 
symmetric case}\label{appt}
Parts (b) and (c) of Lemma~\ref{dontw} for locally symmetric manifolds 
with $\,\rd=1\,$ could, in principle, be derived from Cahen and Parker's 
classification \cite{cahen-parker} of 
pseu\-do-Riem\-ann\-i\-an symmetric manifolds. We prove them here 
directly, for the reader's convenience. Our argument uses assertions (a), 
(d) in Lemma~\ref{dontw}, along with (e) for $\,W\nh$, which were 
established in the proof of Lemma~\ref{dontw} before Appendix II was 
mentioned.

Suppose that $\,\nabla R=0\,$ and $\,\rd=1$. Replacing $\,M\,$ by an open 
subset, we also assume that the \od\ $\,\dz\,$ is spanned by a vector field 
$\,u$. By (\ref{cur}),
\begin{equation}\label{ruo}
\mathrm{i)}\hskip9pt
R(\,\cdot\,,\,\cdot\,)\hs u\,\,=\,\,\varOmega\otimes u\hskip16pt\text{\rm 
or,\ in\ coordinates,}\hskip14pt
\mathrm{ii)}\hskip9ptu^lR_{jkl}{}^s=\varOmega_{jk}u^s,
\end{equation}
for some differential $\,2$-form $\,\varOmega$, which obviously does not 
depend on the choice of $\,u$. (It is also clear from (\ref{cur}) that 
$\,\varOmega\,$ is the curvature form of the connection in the line bundle 
$\,\dz$, induced by the \lcc\ of $\,g$.) Being unique, $\,\varOmega\,$ is 
parallel, and so are $\,\ri\,$ and $\,W\nh$, which implies the Ric\-ci 
identities $\,R\cdot\varOmega=0$, $\,R\cdot\ri=0$, and $\,R\cdot W\nh=0$. In 
coordinates: $\,R_{mlj}{}^s\tau_{sk}+R_{mlk}{}^s\tau_{js}=0$, where 
$\,\tau=\varOmega\,$ or $\,\tau=\ri$, and
\begin{equation}\label{rwp}
R_{qpj}{}^sW_{\nh sklm}\hs+\,R_{qpk}{}^sW_{\nh jslm}
\hs+\,R_{qpl}{}^sW_{\nh jksm}\hs+\,R_{qpk}{}^sW_{\nh jkls}\hs\,=\,\,0\hs.
\end{equation}
Summing $\,R_{mlj}{}^s\varOmega_{sk}+R_{mlk}{}^s\varOmega_{js}=0\,$ 
against $\,u^l\nh$, we obtain $\,\varOmega\circ\hs\varOmega=0$, where the 
metric $\,g\,$ is used to treat $\,\varOmega\,$ as a bundle morphism 
$\,\tm\to\tm\,$ that sends each vector field $\,v\,$ to the vector field 
$\,\varOmega v\,$ with 
$\,g(\varOmega v,v\hh'\hh)=\hs\varOmega(v,v\hh'\hh)\,$ for all vector 
fields $\,v\hh'\nnh$. Lemma~\ref{dontw}(d) and (\ref{ruo}.i) give 
$\,W(\,\cdot\,,\,\cdot\,,u,v)=\hs R(\,\cdot\,,\,\cdot\,,u,v)=0\,$ for our 
fixed vector field $\,u$, spanning $\,\dz$, and any section $\,v\,$ of 
$\,\dzp\nnh$. Hence, by (\ref{wer}), 
$\,g(u,\,\cdot\,)\wedge\sigma(v,\,\cdot\,)
=g(v,\,\cdot\,)\wedge\sigma(u,\,\cdot\,)$. Thus, $\,\sigma u=c\hh u\,$ for 
the Schouten tensor $\,\sigma\,$ and some constant $\,c\hh$, with 
$\,\sigma u\,$ defined analogously to $\,\varOmega v$. (Otherwise, choosing 
$\,v\,$ such that $\,u,\sigma u\,$ and $\,v\,$ are linearly independent at a 
given point $\,x$, we would obtain a contradiction with the equality between 
planes in $\,\txm$, corresponding to the above equality between exterior 
products.) Consequently, 
$\,g(u,\,\cdot\,)\wedge(\sigma+c\hh g)(v,\,\cdot\,)=0$, and so 
$\,\sigma v+c\hh v\,$ is a section of $\,\dz\,$ whenever $\,v\,$ is a section 
of $\,\dzp\nnh$. Let us now fix $\,u\hh'$ as in (\ref{uuo}). Symmetry of 
$\,\sigma\,$ gives $\,g(\sigma u\hh'\nnh,u)=c\hh$. In a suitably ordered 
basis with (\ref{wtt}), at any point $\,x$, the endomorphism of $\,\txm\,$ 
corresponding to $\,\sigma_x$ thus has an upper triangular matrix with the 
diagonal entries $\,c\hh,-\hs c\hh,\dots,-\hs c\hh,c\hh$, so that 
$\,\hs\text{\rm tr}_g\hh\sigma=(4-n)\hs c\hh$. Consequently, 
$\,(n-2)\,\text{\rm s}\hs=2(n-1)(4-n)\hs c\hh$, for the scalar curvature 
$\,\hs\text{\rm s}\hh$, and $\,(n-2)\ri\hs u=2c\hh u$. However, contracting 
(\ref{ruo}.ii) in $\,k=s$, we get $\,\ri\hs u=-\hs\varOmega u$, and so 
$\,(n-2)\hs\varOmega u=-\hs2c\hh u$. The equality 
$\,\varOmega\circ\hs\varOmega=0\,$ that we derived from the Ric\-ci 
identity $\,R\cdot\varOmega=0\,$ now gives $\,c=0$. Hence 
$\,\hs\text{\rm s}\hs=0\,$ (which yields Lemma~\ref{dontw}(c)), and 
$\,\ri\hs u=0$.

As $\,c=0\,$ and $\,\sigma=\ri$, the assertion about $\,\sigma v+c\hh v\,$ 
obtained above means that $\,\ri\hh v\,$ is a section of $\,\dz\,$ 
whenever $\,v\,$ is a section of $\,\dzp\nnh$. Let $\,\lambda,\mu,\xi\,$ 
be the $\,1$-forms with $\,\lambda=g(u,\,\cdot\,)$, 
$\,\mu=g(u\hh'\nnh,\,\cdot\,)$, $\,\xi(u\hh'\hh)=0$, and 
$\,\ri\hh v=\xi(v)\hh u\,$ for sections $\,v\,$ of $\,\dzp\nnh$. 
Transvecting (\ref{ruo}.ii) with $\,\mu_s$, we get 
$\,\varOmega=R(\,\cdot\,,\,\cdot\,,u,u\hh'\hh)
=(n-2)^{-1}\lambda\wedge\ri\hh(u\hh'\nnh,\,\cdot\,)\,$ from 
Lemma~\ref{dontw}(c) with $\,\ri\hs u=0\,$ and Lemma~\ref{dontw}(d). 
However, evaluating $\,\ri\hh(u\hh'\nnh,\,\cdot\,)\,$ on $\,u\hh'\nnh,u\,$ 
and sections $\,v\,$ of $\,\dzp\nnh$, we see that 
$\,\ri\hh(u\hh'\nnh,\,\cdot\,)=h\hs\lambda+\hs\xi$, with 
$\,h=\ri\hh(u\hh'\nnh,u\hh'\hh)$. (Note that $\,\xi(u)=0\,$ since 
$\,\ri\hs u=0$, while $\,\dz\subset\dzp$ by Lemma~\ref{dontw}(a).) 
Therefore,
\begin{equation}\label{omr}
\mathrm{i)}\hskip9pt(n-2)\hs\varOmega\,\,=\,\,\lambda\wedge\hs\xi\hs,
\hskip25pt\mathrm{ii)}\hskip9pt\ri\,\,=\,\,h\hs\lambda\otimes\lambda\,
+\,\lambda\otimes\xi\,+\,\xi\otimes\lambda\hs.
\end{equation}
In addition, if $\,v\hh'$ denotes the unique vector field with 
$\,g(v\hh'\nnh,\,\cdot\,)=\xi$, then $\,u\,$ and $\,v\hh'$ are null and 
orthogonal, or, equivalently,
\begin{equation}\label{nul}
\text{\rm the\ $\,1$-forms\ $\,\lambda\,$\ and\ $\,\xi\,$\ are\ null\ and\ 
mutually\ orthogonal.}
\end{equation}
In fact, $\,g(u,u)=0\,$ by Lemma~\ref{dontw}(a), 
$\,g(u,v\hh'\hh)=0\,$ as $\,\xi(u)=0$, and $\,v\hh'$ is null since 
(\ref{omr}) yields $\,(n-2)\hs[\ri\hs(\varOmega\hs u\hh'\hh)
-\varOmega\hs(\ri\hs u\hh'\hh)]=2\hh g(v\hh'\nnh,v\hh'\hh)\hs u$, 
while, transvecting the Ric\-ci identity 
$\,R_{mlj}{}^sR_{sk}+R_{mlk}{}^sR_{js}=0\,$ with $\,u^l$ and using 
(\ref{ruo}.ii), we see that $\,\ri\,$ and $\,\varOmega$ commute 
as bundle morphisms $\,\tm\to\tm$.

Furthermore, transvecting with $\,\mu^k\mu^m$ the coordinate form 
$\,R_{mlj}{}^s\tau_{sk}+R_{mlk}{}^s\tau_{js}=0$ of the Ric\-ci 
identity $\,R\cdot\tau=0\,$ for the parallel tensor field 
$\,\tau=(n-2)\hs\varOmega+\ri
=h\hs\lambda\otimes\lambda+2\lambda\otimes\xi\,$ (cf.\ (\ref{omr})), 
we get $\,2\hh\lambda_jb_{\hs ls}\xi^s=0$, where 
$\,b=W(u\hh'\nnh,\,\cdot\,,u\hh'\nnh,\,\cdot\,)$. Namely, 
$\,R=\hs W\hs+\,(n-2)^{-1}\hs g\wedge\ri\,$ by Lemma~\ref{dontw}(c), 
$\,W_{\nh mlj}{}^s\tau_{sk}=0\,$ in view of Lemma~\ref{dontw}(d), 
$\,\mu^k\mu^mW_{\nh mlk}{}^s\tau_{js}=2\hh\lambda_jb_{\hs ls}\xi^s$ since 
$\,b\hs(u,\,\cdot\,)=0\,$ (again from Lemma~\ref{dontw}(d)), and the 
remaining terms, related to $\,g\wedge\ri$, add up to $\,0\,$ as a 
consequence of (\ref{nul}), (\ref{omr}.ii) and the formula for $\,\tau$. 
(Note that (\ref{nul}) gives $\,R_j{}^s\tau_{sk}=R_j{}^s\tau_{ks}=0$, 
and so four out of the eight remaining terms vanish individually.) 
However, $\,u\ne0$, and so $\,\lambda\ne0$, which gives 
$\,b\hs(\,\cdot\,,v\hh'\hh)=0$, where $\,v\hh'$ is the vector field 
with $\,g(v\hh'\nnh,\,\cdot\,)=\xi$. Thus, 
$\,W(u\hh'\nnh,\,\cdot\,,u\hh'\nnh,v\hh'\hh)=0$. As a result, the 
$\,3$-ten\-sor $\,W(\,\cdot\,,\,\cdot\,,\,\cdot\,,v\hh'\hh)\,$ must vanish: 
it yields the value $\,0\,$ whenever each of the three arguments is 
either $\,u\hh'$ or a section of $\,\dzp\nnh$. (Lemma~\ref{dontw}(e) for 
$\,W\hs$ is already established.)

The relation $\,W(\,\cdot\,,\,\cdot\,,\,\cdot\,,v\hh'\hh)=0\,$ implies in 
turn that $\,W(\,\cdot\,,\,\cdot\,,\,\cdot\,,\ri\hs v)=0\,$ (in 
coordinates: $\,W_{\nh jkl}{}^sR_{sp}=0$). In fact, by (\ref{omr}.ii), the 
image of $\,\ri\,$ is spanned by $\,u\,$ and $\,v\hh'\nnh$, while 
$\,W(\,\cdot\,,\,\cdot\,,\,\cdot\,,u)=0\,$ according to 
Lemma~\ref{dontw}(d).

As in \cite[1$^{\hs\text{\rm o}}$ on p.~214]{olszak}, we have 
$\,W\nh=(\lambda\otimes\lambda)\wedge\hs b\,$ (notation of (\ref{wer})), 
where, again, $\,b=W(u\hh'\nnh,\,\cdot\,,u\hh'\nnh,\,\cdot\,)$. Namely, by 
Lemma~\ref{dontw}(e) for $\,W\nh$, both sides agree on any quadruple of vector 
fields, each of which is either $\,u\hh'$ or a section of $\,\dzp\nnh$.

Finally, transvecting (\ref{rwp}) with $\,\mu^k\mu^m$ and replacing 
$\,R\,$ by $\,W\nh+(n-2)^{-1}\hs g\wedge\ri$, we obtain two 
contributions, one from $\,W\hs$ and one from $\,g\wedge\ri$, the sum 
of which is zero. Since $\,W\nh=(\lambda\otimes\lambda)\wedge\hs b$, the 
$\,W\hs$ contribution vanishes: its first two terms add up to $\,0$, and so do 
its other two terms. (As we saw, $\,b\hs(u,\,\cdot\,)=0$, while, obviously, 
$\,b\hs(u\hh'\nnh,\,\cdot\,)=0$.) Out of the sixteen terms forming the 
$\,g\wedge\ri\,$ contribution, eight are separately equal to zero since  
$\,W_{\nh jkl}{}^sR_{sp}=0$, and so, in view of (\ref{omr}.ii) and the 
relation  $\,W\nh=(\lambda\otimes\lambda)\wedge\hs b\hh$, vanishing of the 
$\,g\wedge\ri\,$ contribution gives $\,\lambda_pS_{jlq}=\lambda_qS_{jlp}$, for 
$\,S_{jlq}=2\hh b_{jl}\xi_q-\hh b_{ql}\xi_j
-\hh b_{qj}\xi_l$. Thus, $\,S_{jlq}=\eta_{jl}\lambda_q$ for some 
twice-co\-var\-i\-ant symmetric tensor field $\,\eta$, which, summed 
cyclically over $\,j,l,q$, yields $\,0\,$ (due to the definition 
of $\,S_{jlq}$ and symmetry of $\,b$). As $\,\lambda\ne0\,$ and the 
symmetric product has no zero divisors, we get $\,\eta=0\,$ and 
$\,S_{jlq}=0$. The expression $\,b_{jl}\xi_q-\hh b_{ql}\xi_j$ is, 
therefore, skew-sym\-met\-ric in $\,j,l$. As it is also, clearly, 
skew-sym\-met\-ric in $\,j,q$, it must be totally skew-sym\-met\-ric 
and hence equal to one-third of its cyclic sum over $\,j,l,q$. That 
cyclic sum, however, is $\,0\,$ in view of symmetry of $\,b$, so that 
$\,b_{jl}\xi_q=\hh b_{ql}\xi_j$. Thus, $\,\xi=0$, for otherwise the last 
equality would yield $\,b\hs=\varphi\hs\xi\otimes\hs\xi\,$ for some function 
$\,\varphi$, and hence $\,W\nh=(\lambda\otimes\lambda)\wedge\hs b\hs
=\hh\varphi\hs(\lambda\otimes\lambda)\wedge(\hh\xi\otimes\hs\xi)$, which would 
clearly imply that the vector field $\,v\hh'$ with 
$\,g(v\hh'\nnh,\,\cdot\,)=\xi\,$ is a section of the \od\ $\,\dz$, not 
equal to a function times $\,u\,$ (as $\,\xi(u\hh'\hh)=0$, while 
$\,g(u,u\hh'\hh)=1$), contradicting one-di\-men\-sion\-al\-i\-ty of $\,\dz$. 
Therefore, $\,\ri=h\hs\lambda\otimes\lambda\,$ by (\ref{omr}.ii) with 
$\,\xi=0$, which proves assertion (b) of Lemma~\ref{dontw} in our case.


\begin{thebibliography}{99} 


\bibitem{cahen-parker}\textsc{Cahen, M. \&\ Parker, M.}, 
Pseu\-do-riem\-ann\-i\-an symmetric spaces. {\em Mem.\ Amer.\ Math.\ 
Soc.} \textbf{229} (1980), 1--108.

\bibitem{chaki-gupta}\textsc{Chaki, M.\hskip1.8ptC. \&\ Gupta, B.}, On 
con\-for\-mal\-ly symmetric spaces. {\em Indian J.\ Math.} \textbf{5} (1963), 
\hbox{113\hs--}122.

\bibitem{derdzinski-78}\textsc{Derdzi\'nski, A.}, On homogeneous 
con\-for\-mal\-ly symmetric pseudo-Riemannian manifolds. {\em Colloq.\ Math.} 
\textbf{40} (1978), 167--185.

\bibitem{derdzinski-roter-77}\textsc{Derdzi\'nski, A. \&\ Roter, W.}, On 
con\-for\-mal\-ly symmetric manifolds with metrics of indices $\,0\,$ and 
$\,1\hh$. {\em Tensor\/} (N.\hskip1.9ptS.) \textbf{31} (1977), 255\hs--259.

\bibitem{derdzinski-roter-78}\textsc{Derdzi\'nski, A. \&\ Roter, W.}, 
Some theorems on con\-for\-mal\-ly symmetric manifolds. {\em Tensor\/} 
(N.\hskip1.9ptS.) \textbf{32} (1978), 11--23.

\bibitem{derdzinski-roter-80}\textsc{Derdzi\'nski, A. \&\ Roter, W.}, 
Some properties of con\-for\-mal\-ly symmetric manifolds which are not 
Ric\-ci-re\-cur\-rent. {\em Tensor\/} (N.\hskip1.9ptS.) \textbf{34} (1980), 
11--20.

\bibitem{derdzinski-roter-07}\textsc{Derdzinski, A. \&\ Roter, W.}, 
Projectively flat surfaces, null parallel distributions, and 
con\-for\-mal\-ly symmetric manifolds. Preprint, math.DG/0604568. To appear in 
{\em Tohoku Math.\ J}. 

\bibitem{derdzinski-roter}\textsc{Derdzinski, A. \&\ Roter, W.}, 
Compact pseu\-do\hs-Riem\-ann\-i\-an manifolds with parallel Weyl 
tensor. Preprint, http:/\hskip-1.5pt/arXiv.org/abs/math.DG/0702491.

\bibitem{deszcz}\textsc{Deszcz, R.}, On hypercylinders in conformally 
symmetric manifolds. {\em Publ.\ Inst.\ Math.} (Beograd) 
(N.\hskip1.9ptS.) \textbf{51}(65) (1992), 101--114.

\bibitem{deszcz-hotlos}\textsc{Deszcz, R. \&\ Hotlo\'s, M.}, On a certain 
subclass of pseudosymmetric manifolds. {\em Publ.\ Math.\ Debrecen\/} 
\textbf{53} (1998), 29\hh--48.

\bibitem{dillen-verstraelen}\textsc{Dillen, F.\hskip1.8ptJ.\hskip1.8pt E. \&\ 
Verstraelen, L.\hskip1.8ptC.\hskip1.8ptA.} (eds.), {\em Handbook of 
Differential Geometry, Vol.~I}. North-Hol\-land, Amsterdam, 2000.

\bibitem{hotlos}\textsc{Hotlo\'s, M.}, On conformally symmetric warped 
products. {\em Ann.\ Acad.\ Paedagog.\ Cracov.\ Stud.\ Math.} \textbf{4} 
(2004), 75\hs--85.

\bibitem{olszak}\textsc{Olszak, Z.}, On con\-for\-mal\-ly recurrent 
manifolds, I: Special distributions. {\em Zesz. Nauk. Po\-li\-tech. \'Sl., 
Mat.-Fiz.} \textbf{68} (1993), 213\hs--225.

\bibitem{patterson-walker}\textsc{Patterson, E.\hskip1.8ptM. \&\ Walker, 
A.\hskip1.6ptG.}, Riemann extensions. {\em Quart.\ J.\ Math.\ Oxford\/} Ser.\ 
(2) \textbf{3} (1952), 19\hs--28.

\bibitem{rong}\textsc{Rong, J.\hskip1.8ptP.}, On $^2\nh K^*_n$ space. {\em 
Tensor\/} (N.\hskip1.9ptS.) \textbf{49} (1990), 117--123.

\bibitem{roter}\textsc{Roter, W.}, On con\-for\-mal\-ly symmetric 
Ric\-ci-re\-cur\-rent spaces. {\em Colloq.\ Math.} \textbf{31} (1974), 
87--\hh96.

\bibitem{sharma}\textsc{Sharma, R.}, Proper conformal symmetries of conformal 
symmetric space-times. {\em J.\ Math.\ Phys.} \textbf{29} (1988), 2421--2422.

\bibitem{simon}\textsc{Simon, U.}, Compact con\-for\-mal\-ly symmetric 
Riemannian spaces. {\em Math.\ Z.} \textbf{132} (1973), \hbox{173\hs--}177.

\end{thebibliography}
\end{document}